\documentclass[12pt]{amsart}

\title{$L$-smooth factorization for Noetherian $F$-finite rings}
\author{Manuel Blickle and Daniel Fink}

\usepackage[margin=2.5cm]{geometry}
\usepackage{hyperref}
\usepackage{times}
\usepackage[utf8]{inputenc}
\usepackage[T1]{fontenc}
\usepackage{amsmath, amssymb, amsthm, mathtools}

\usepackage{tikz}
\usetikzlibrary{cd,babel,graphs}

\theoremstyle{definition}
\newtheorem{defi}{Definition}[section]
\newtheorem{ex}[defi]{Example}
\newtheorem{con}[defi]{Construction}
\newtheorem{rem}[defi]{Remark}

\theoremstyle{plain}
\newtheorem{thm}[defi]{Theorem}
\newtheorem{prop}[defi]{Proposition}
\newtheorem{lem}[defi]{Lemma}
\newtheorem{cor}[defi]{Corollary}

\newcommand{\NN}{\mathbb{N}}  

\newcommand\FF{\mathbb{F}}

\newcommand{\spec}{\operatorname{Spec}}

\newcommand{\Tor}{\operatorname{Tor}}

\newcommand{\id}{\operatorname{id}}

\newcommand{\comp}{\operatorname{Comp}}
\newcommand{\sym}{\operatorname{Sym}}
\newcommand{\gr}{\operatorname{gr}}

\renewcommand{\to}[1][]{\xrightarrow{\ #1\ }}

\newcommand{\into}[1][]{\xhookrightarrow{\ #1\ }}

\usepackage[backend=biber,style=alphabetic, natbib=true,hyperref=true,language=auto, autolang=other,maxnames=50, sortcase=false, minalphanames=4, maxalphanames=4]{biblatex}
\nocite{*}

\begin{document}

\begin{abstract}
We show that any homomorphism between Noetherian $F$-finite rings can be factored into a regular morphism between Noetherian $F$-finite rings followed by a surjection. This result establishes an analog of the ‘smooth-by-surjective’ factorization for finite type maps. As part of our analysis, we observe that for maps of Noetherian $F$-finite rings, regularity and formal smoothness are both equivalent to \emph{$L$-smoothness}, meaning that the cotangent complex, as in the smooth case, is a locally free module of finite rank concentrated in degree zero.

Our findings may also be viewed as a relative version of Gabber’s final remark in \citep{Gab04}, which states that any Noetherian $F$-finite ring is a quotient of a regular Noetherian $F$-finite ring.
\end{abstract}

\maketitle


\section{Introduction}
Noetherian $\mathbb{F}_p$-algebras that are $F$-finite, meaning their Frobenius endomorphism is finite, share many properties with algebras of finite type over a field. For instance, these rings have finite Krull dimension \citep{Kun76}, are excellent \citep{Kun76}, can be expressed as quotients of regular Noetherian $F$-finite $\FF_p$-algebras \citep{Gab04}, and possess dualizing complexes \citep{Gab04,BBST}.

The defining property of finite type maps of rings ensures that any morphism of this kind can be factored into a smooth map followed by a surjective map. This is fundamental, as it often reduces problems about finite type morphisms to those about smooth maps and closed immersions. The aim of this article is to establish an analogous factorization for arbitrary maps of Noetherian $F$-finite rings. Since such rings are ubiquitous in the study of commutative algebra and algebraic geometry in positive characteristic, this result applies to a broad class of morphisms. For instance, any perfect field is Noetherian and $F$-finite, and these properties extend to algebras that are essentially of finite type over a Noetherian $F$-finite ring, as well as to completions of Noetherian $F$-finite rings.

Since Noetherian $F$-finite rings and their morphisms are generally not of finite type, we cannot rely on the usual finiteness conditions central to the definition of smoothness. Instead, we adopt a more flexible homological notion inspired by a result of André and Quillen: a map of finite presentation is smooth if and only if its cotangent complex is a finitely generated projective module concentrated in degree zero \citep[]{Qui70}. Motivated by this characterization, we define a ring homomorphism to be \emph{$L$-smooth} if its cotangent complex resembles that of a smooth map; see Definition \ref{lsmoothdefi}. This notion proves particularly useful in the context of this article because it aligns closely with an intrinsic property of Noetherian $F$-finite rings: these algebras can be characterized by their “absolute” module of Kähler differentials. Specifically, a Noetherian $\FF_p$-algebra $R$ is $F$-finite if and only if $\Omega_{R/\FF_p}$ is a finitely generated $R$-module \citep{Fog80}; see Section 2 for further discussion.

The main result of this paper is the following Theorem; see Theorem \ref{maintheorem} in the main text.

\begin{thm}
\label{introthm}
Any map $R\to S$ of Noetherian $F$-finite $\FF_p$-algebras admits a factorization
\[
\begin{tikzcd}
                         & {T} \arrow[rd, "\psi"] &    \\
{R} \arrow[ru, "\varphi"] \arrow[rr] &               & {S}
\end{tikzcd}
\]
in which $T$ is a Noetherian $F$-finite $\FF_p$-algebra, $\varphi$ is $L$-smooth, and $\psi$ is surjective.
\end{thm}

\subsection*{Notation and Conventions} Throughout this paper, all rings and algebras are assumed to be commutative, associative and unital, and all ring homomorphisms are assumed to preserve the identity.

For a ring $R$, we denote by $D(R)$ the derived $\infty$-category of $R$-modules, as discussed in \citep{Lur17}. However, except in the final chapter, the article does not require prior knowledge of $\infty$-categories. Whenever inverse limits appear, they are to be interpreted as derived limits, as described in \cite[\href{https://stacks.math.columbia.edu/tag/08TB}{Section 08TB}]{stacks-project}. Moreover, fiber sequences should be understood as distinguished triangles. 

Cohomological indexing is used throughout.

For an $\FF_p$-algebra, the Frobenius homomorphism is denoted by $F$. 

Given a ring map $f: R \to S$, the associated cotangent complex is typically denoted $L_{S/R}$. To avoid ambiguity, however, the notation $L_f$ is sometimes used, especially when referring to the cotangent complex of the Frobenius homomorphism.

\subsection*{Outline} In Section 2, we collect basic results for $L$-smooth maps and compare this notion with familiar properties such as smoothness, formal smoothness, and regularity. We then focus on maps between Noetherian $F$-finite rings, where we show that the notions of $L$-smoothness, formal smoothness, and regularity coincide.

In Section 3, we show the existence of a factorization as described in Theorem \ref{introthm} for maps of \emph{regular} Noetherian $F$-finite rings; see Corollary \ref{regfact}.

In Section 4, we build on the results of Section 3 to prove Theorem \ref{introthm} in its full generality; see Theorem \ref{maintheorem}.

Finally, in Section 5, we provide an alternative perspective on the factorization using the theory of Adams completion as developed in \citep{BBST}.

\subsection*{Acknowledgements} 
The problem discussed in this article originated from a conversation between the first author and Bhargav Bhatt, who suggested a solution along the lines of the final chapter, drawing on ideas from \cite{BBST}. At the time, it was unclear whether a more classical approach could address this factorization result. Both authors are deeply grateful to Bhatt for his inspiring suggestions and generosity in sharing his ideas.

The second author thanks Eamon Quinlan-Gallego for valuable discussions on this topic. He also thanks the authors of \cite{BBST} for sharing an early manuscript of their work, which provided significant insights.

The authors acknowledge support from the Deutsche Forschungsgemeinschaft (DFG, German Research Foundation) through the Collaborative Research Centre TRR 326 \textit{Geometry and Arithmetic of Uniformized Structures}, project number 444845124.

\section{L-smoothness of ring maps}
Understanding the relationship between \emph{$L$-smoothness} and the classical notions, such as smoothness, formal smoothness, and regularity, is a central goal of this chapter. To this end, we first recall that a homomorphism of rings $R\to S$ is called \emph{formally smooth}, if for any commutative diagram 
\[
\begin{tikzcd}
{R} \arrow[r] \arrow[d] & {A} \arrow[d] \\
{S} \arrow[r]           & {A/I}          
\end{tikzcd}
\]
where $I\subset A$ is an ideal with $I^2=0$, there exists a morphism $S\to A$ that makes the entire diagram commute. Furthermore, if $R\to S$ is also of finite presentation, it is called \emph{smooth}. 

André and Quillen \citep[Theorems 5.4 and 5.5]{Qui70} provided a homological characterization of smoothness using the cotangent complex. For completeness, we include a proof of this well-known result, as the details are omitted in the cited reference. 

\begin{thm}
\label{motivation}
Let $R\to S$ be a ring map of finite presentation. The following statements are equivalent:
\begin{enumerate}
\item $R\to S$ is smooth.
\item $L_{S/R}\simeq \Omega_{S/R}[0]\in D(S)$ and $\Omega_{S/R}$ is a finitely generated projective $S$-module.
\end{enumerate}
\end{thm}
\begin{proof}
Since $R\to S$ is of finite presentation, it follows that $\Omega_{S/R}$ is a finitely presented $S$-module; see \cite[\href{https://stacks.math.columbia.edu/tag/00RY}{Lemma 00RY}]{stacks-project}.

First, assume that $R\to S$ is smooth. By \cite[\href{https://stacks.math.columbia.edu/tag/08R5}{Lemma 08R5}]{stacks-project}, we have $L_{S/R}\simeq\Omega_{S/R}[0]$. Furthermore, $\Omega_{S/R}$ is projective by \cite[\href{https://stacks.math.columbia.edu/tag/031J}{Proposition 031J}]{stacks-project}. 

Conversely, suppose that $L_{S/R}\simeq \Omega_{S/R}[0]$ and that $\Omega_{S/R}$ is a projective $S$-module. By \cite[\href{https://stacks.math.columbia.edu/tag/08RB}{Lemma 08RB}]{stacks-project}, it follows that the cotangent complex and the naive cotangent complex are equivalent. Consequently, $R\to S$ is formally smooth by \cite[\href{https://stacks.math.columbia.edu/tag/031J}{Proposition 031J}]{stacks-project}. Since $R\to S$ is also finitely presented by assumption, it follows that $R\to S$ is smooth.
\end{proof}

This motivates our definition of $L$-smoothness:

\begin{defi}
\label{lsmoothdefi}
A ring homomorphism $R\to S$ is called \emph{$L$-smooth}, if its cotangent complex $L_{S/R}$ is concentrated in homological degree zero, and $H^0(L_{S/R})\cong \Omega_{S/R}$ is a finitely generated projective $S$-module.
\end{defi}

The proof of Theorem \ref{motivation} shows the following relations:

\begin{lem}
For a ring map $f: R \to S$, the following implications hold:
\[
    f\text{ is smooth} \implies f\text{ is $L$-smooth} \implies f\text{ is formally smooth}
\]
\end{lem}
The reverse of each implication is not true as the following two examples show.
\begin{ex}
\label{formallysmooth} 
Let $R$ be any ring and let $S$ be a polynomial ring in infinitely many variables over $R$. Then the canonical map $R\to S$ is formally smooth, but $\Omega_{S/R}$ is not finitely generated as an $S$-module. Hence, the map is not $L$-smooth.
\end{ex}

\begin{ex}
Let $k$ be a perfect ring of infinite cardinality. Since $k$ is perfect, the cotangent complex $L_{k/\FF_p}$ vanishes by \cite[\href{https://stacks.math.columbia.edu/tag/0G60}{Lemma 0G60}]{stacks-project}. Consequently, the canonical map $\FF_p\to k$ is $L$-smooth. However, this map is not of finite type and therefore not smooth.
\end{ex}

Let $R$ be a ring and $k\ge0$ a natural number. An object $C\in D(R)^{\le0}$ is said to have \emph{flat dimension $\leq k$} if it is equivalent to a complex $Q^{\bullet}$ where each $Q^i$ is flat and $Q^i=0$ for $i\notin [-k,0]$. The flat dimension of $C$ is the smallest $k$ such that $C$ has flat dimension $\leq k$. 

Recall that a right-bounded complex is called \emph{pseudo-coherent} if it is equivalent to a right-bounded complex consisting of finitely generated free modules
\cite[\href{https://stacks.math.columbia.edu/tag/064Q}{Definition 064Q}]{stacks-project}.
The following lemma provides a basic but useful criterion for $L$-smoothness:

\begin{lem}
\label{equivlsmooth}
A morphism $R\to S$ of rings is $L$-smooth if and only if its cotangent complex $L_{S/R}$ is pseudo-coherent with flat dimension $0$.
\end{lem}
\begin{proof}
This result is shown in \cite[\href{https://stacks.math.columbia.edu/tag/0658}{Lemma 0658}]{stacks-project}.
\end{proof}

\begin{prop}
\label{compositionlsmooth}
Let $R\to S$ and $S\to T$ be $L$-smooth ring maps. Then the composition $R\to T$ is also $L$-smooth.
\end{prop}
\begin{proof}
This follows from the fundamental fiber sequence $L_{S/R}\otimes^L_S T\to L_{T/R} \to L_{T/S}$.
Indeed, by \cite[\href{https://stacks.math.columbia.edu/tag/064R}{Lemma 064R}]{stacks-project} $L_{T/R}$ is pseudo coherent. Furthermore, by \cite[\href{https://stacks.math.columbia.edu/tag/0655}{Lemma 0655}]{stacks-project} it has flat dimension $0$.
\end{proof}

\begin{prop}
\label{basechangelsmooth}
Let
\[
\begin{tikzcd}
{R'} \arrow[r] \arrow[d] & {S'} \arrow[d] \\
{R} \arrow[r]           & {S}          
\end{tikzcd}
\]
be a pushout square of commutative rings, and suppose that $R'\to S'$ is $L$-smooth. If $R$ and $S'$ are Tor-independent over $R'$, then $R\to S$ is also $L$-smooth. 
\end{prop}
\begin{proof}
The Tor-independence of $R$ and $S'$ over $R'$ implies that the diagram above is a pushout square of animated rings. Consequently, by \citep[Remark 25.3.2.4]{Lur18}, we obtain a natural equivalence
\[
L_{R/S}\simeq L_{R'/S'}\otimes^L_{S'} S
\]
in $D(S)$. Alternatively, see \cite[\href{https://stacks.math.columbia.edu/tag/08QQ}{Lemma 08QQ}]{stacks-project}. Now, $L_{R/S}$ is pseudo-coherent, and by \cite[\href{https://stacks.math.columbia.edu/tag/066L}{Lemma 066L}]{stacks-project}, its flat dimension is $0$.
\end{proof}

\begin{rem}
Since $L$-smooth morphisms are formally smooth, it follows from \citep[Théorème 19.7.1]{EGAIV}, that $L$-smooth morphisms of Noetherian rings are flat. Hence, if $R'$ and $S'$ in the previous proposition are Noetherian, then $R$ and $S'$ are automatically Tor-independent over $R'$.
\end{rem}

The cotangent complex can also be used to classify other properties of ring maps when the rings involved are Noetherian. Recall that a map of Noetherian rings $R\to S$ is called \emph{regular} if it is flat and, for every prime ideal $p\in \spec(R)$, and for every finite purely inseparable field extension $k/\kappa(p)$, the ring $S\otimes_R k$ is regular; see \cite[\href{https://stacks.math.columbia.edu/tag/07BZ}{Definition 07BZ}]{stacks-project}.

Due to André \citep[Supplement Theorem 30]{And74}, we have the following characterization of regularity.

\begin{thm}
Let $R\to S$ be a morphism of Noetherian rings. The following statements are equivalent:
\begin{enumerate}
\item $R\to S$ is regular.
\item The cotangent complex $L_{S/R}$ is equivalent to $\Omega_{S/R}[0]\in D(S)$, and $\Omega_{S/R}$ is a flat $S$-module.
\end{enumerate}
\end{thm}

We also recall the following characterization of formal smoothness, which fundamentally builds on the contributions of Avramov, Briggs, and Iyengar regarding the cotangent complex.

\begin{thm}
Let $R\to S$ be a map of Noetherian rings. The following statements are equivalent:
\begin{enumerate}
\item $R\to S$ is formally smooth.
\item The cotangent complex $L_{S/R}$ is equivalent to $\Omega_{S/R}[0]\in D(S)$, and $\Omega_{S/R}$ is a projective $S$-module.
\end{enumerate}
\end{thm}
\begin{proof}
The implication (2) $\implies$ (1) was shown in the proof of Theorem \ref{motivation} and does not require the Noetherian hypothesis. 

Conversely, suppose that $R\to S$ is formally smooth. By \cite[\href{https://stacks.math.columbia.edu/tag/031J}{Proposition 031J}]{stacks-project}, part (6), this implies that $\Omega_{S/R}$ is projective over $S$. Next, we show that $L_{S/R}$ is concentrated in homological degree $0$. 
Consider a factorization $R \to P \to S$, where $P$ is a polynomial ring over $R$, and $P\to S$ is surjective. By \cite[\href{https://stacks.math.columbia.edu/tag/08QH}{Lemma 08QH}]{stacks-project}, we have $L_{P/R}\simeq \Omega_{P/R}[0]$. Since $\Omega_{P/R}$ is a free $P$-module, it follows that
\[
L_{P/R}\otimes_P^L S\simeq (\Omega_{P/R}\otimes_P S)[0].
\]
By \cite[\href{https://stacks.math.columbia.edu/tag/08RA}{Lemma 08RA}]{stacks-project}, the first cohomology groups of $L_{S/P}$ satisfy $H^0(L_{S/P})=0$ and $H^{-1}(L_{S/P})=I/I^2$, where $I$ is the kernel of the surjection $P\to S$. Moreover, by \cite[\href{https://stacks.math.columbia.edu/tag/031J}{Proposition 031J}]{stacks-project}, part (5), we have a split exact sequence:
\[
0 \to H^{-1}(L_{S/P}) \to H^0(L_{P/R}\otimes_P^L S) \to H^0(L_{S/R}) \to 0
\]
From these observations, we deduce that, for any $S$-module $M$, the canonical fiber sequence in $D(S)$,
\[
L_{P/R}\otimes_P^L S\otimes_S^L M \to L_{S/R}\otimes_S^L M \to L_{S/P}\otimes_S^L M,
\]
induces the following exact sequence of $S$-modules:
\[
0\to H^{-1}(L_{S/R}\otimes_S^L M)\to H^{-1}(L_{S/P})\otimes_S M \to H^0(L_{P/R}\otimes_P^L S)\otimes_S M \to H^0(L_{S/R})\otimes_S M \to 0
\]
Since the map $H^{-1}(L_{S/P}) \to H^0(L_{P/R}\otimes_P^L S)$ is a split injection, it remains injective after tensoring with $M$. As a result, we find that $H^{-1}(L_{S/R}\otimes_S^L M)=0$. Therefore, by \citep[Theorem A]{BI23}, $R\to S$ is a local complete intersection morphism. By \citep[(1.2) Second vanishing theorem]{Avr99}, this implies that the flat dimension of $L_{S/R}$ is at most $1$. In particular, $L_{S/R}$ is concentrated in degrees $-1$ and $0$. However, we have already seen that $H^{-1}(L_{S/R})=0$, so it follows that $L_{S/R}\simeq \Omega_{S/R}[0]$.
\end{proof}

In particular, for maps between Noetherian rings, formal smoothness implies regularity. Moreover, for classes of ring maps where the Kähler differentials are finitely generated, formal smoothness and regularity are equivalent to $L$-smoothness. This equivalence holds, for example, for all maps of finite type between Noetherian rings. In the following, we describe another significant class of maps for which this equivalence persists. 

\subsection*{F-finite ring maps}
Throughout the remainder of this section, all rings are assumed to be $\FF_p$-algebras.

Let $f: R \to S$ be a morphism of $\FF_p$-algebras. The commutative diagram
\[
\begin{tikzcd}
{R} \arrow[r, "F"] \arrow[d, "f"] & {R} \arrow[d, "f"] \\
{S} \arrow[r, "F"] & {S}          
\end{tikzcd}
\]
induces a unique ring homomorphism
\[
F_{S/R}: S \otimes_R F_* R \to S,
\]
called the \emph{relative Frobenius homomorphism} of $f$. Here, $F_*$ denotes the restriction of scalars along the (absolute) Frobenius $F: R \to R $. 

Following \citep{Has15}, we define $F$-finite morphisms using the relative Frobenius homomorphism:
\begin{defi}
A morphism $f: R \to S$ is called \emph{$F$-finite} if its relative Frobenius homomorphism is a finite ring map.
\end{defi}
Equivalently, $f$ is $F$-finite if $S$ is a finite module over its subring $F_{S/R}(S\otimes_R F_*R)=S^p[f(R)]$.


\begin{ex}
An $\FF_p$-algebra $R$ is $F$-finite, if and only if its structure morphism $\FF_p\to R$ is $F$-finite. Moreover, any map to an $F$-finite $\FF_p$-algebra is $F$-finite; see \citep[Lemma 2]{Has15} for these and other basic properties.
\end{ex}

\begin{rem}
Let  $R \to S$  be a morphism of $\FF_p$-algebras, and consider the exact sequence
\[
\Omega_{R[S^p]/R}\otimes_{R[S^p]}S\to \Omega_{S/R}\to \Omega_{S/R[S^p]}\to 0
\]
associated with $R \to R[S^p] \into S$. Since the first map in this sequence is zero, it induces an isomorphism $\Omega_{S/R[S^p]} \cong \Omega_{S/R}$. Hence, $\Omega_{S/R}$ is a finitely generated $S$-module whenever $R\to S$ is $F$-finite.

For maps of Noetherian rings, the converse also holds: a map $R\to S$ is $F$-finite if and only if $\Omega_{S/R}$ is a finitely generated $S$-module; see \citep[Proposition 1]{Fog80} and \citep[Remark 13]{Has15}.
\end{rem}

Thus, we obtain the following:

\begin{prop}
For an $F$-finite map of Noetherian rings $R\to S$, both regularity and formal smoothness are equivalent to $L$-smoothness. Furthermore, any $L$-smooth map of Noetherian $\FF_p$-algebras is $F$-finite.
\end{prop}

\begin{rem}
With this result at hand, one could define a map of Noetherian $\FF_p$-algebras to be \emph{$F$-smooth} if it is $F$-finite and formally smooth. However, since this notion is equivalent to $L$-smoothness, we do not adopt it here.
\end{rem}

\begin{cor}
\label{regularislsmooth}
A Noetherian $F$-finite $\FF_p$-algebra $R$ is regular if and only if it is $L$-smooth over $\FF_p$.
\end{cor}
\begin{proof}
Since $\FF_p$ is a perfect field, $R$ is regular if and only if $\FF_p\to R$ is a regular ring map.
\end{proof}

Cotangent complexes of maps between $F$-finite rings possess a strong finiteness property (see also \citep[Remark 3.3.3]{BBST}, where this is stated for the absolute cotangent complex). The following result generalizes the fact that modules of Kähler differentials for maps between $F$-finite rings are finitely generated. Its proof follows a similar approach.

\begin{lem}
\label{pseudocoh}
Let $R\to S$ be a morphism of Noetherian $F$-finite $\FF_p$-algebras. Then its cotangent complex $L_{S/R}$ is pseudo-coherent.
\end{lem}
\begin{proof}
The sequence of maps $\FF_p\to R\to S$ induces the fiber sequence 
\[
L_{R/\FF_p}\otimes^L_{R}S\to L_{S/\FF_p}\to L_{S/R}.
\]
By \cite[\href{https://stacks.math.columbia.edu/tag/0650}{Lemma 0650}]{stacks-project} and the “two-out-of-three” property of pseudo-coherence in a fiber sequence \cite[\href{https://stacks.math.columbia.edu/tag/064V}{Lemma 064V}]{stacks-project}, it suffices to show that $L_{R/\FF_p}$ is pseudo-coherent if $R$ is an $F$-finite $\FF_p$-algebra. For this, we consider the transitivity fiber sequence
\[
L_{R/\FF_p}\otimes^L_{R}F_*R\to L_{R/\FF_p}\to L_{F}
\]
induced by $\FF_p\to R\to[F] R$. Since the first map in this fiber sequence is a zero map, we obtain a canonical equivalence
\[
L_F\simeq (L_{R/\FF_p}\otimes^L_{R}F_*R)[1]\oplus L_{R/\FF_p}.
\]
Consequently, by \cite[\href{https://stacks.math.columbia.edu/tag/064X}{Lemma 064X}]{stacks-project}, it suffices to show that $L_F$ is pseudo-coherent. But this is true by \citep[Proposition 4.12]{Qui70}, as $R$ is Noetherian and $F$-finite by assumption.
\end{proof}


Using Fogarty's result, we can characterize $F$-finiteness of Noetherian rings in terms of their absolute cotangent complex:

\begin{prop}
Let $R$ be a Noetherian $\FF_p$-algebra. Then $R$ is $F$-finite if and only if $L_{R/\FF_p}$ is pseudo-coherent.
\end{prop}
\begin{proof}
The “only if” direction is a special case of the previous lemma. Conversely, if $L_{R/\FF_p}$ is pseudo-coherent, then by \cite[\href{https://stacks.math.columbia.edu/tag/066E}{Lemma 066E}]{stacks-project}, $H^0(L_{R/\FF_p})=\Omega_{R/\FF_p}$ is a finite $R$-module. Therefore, $R$ is $F$-finite by \citep[Proposition 1]{Fog80}.
\end{proof}

\begin{rem}
Let $R$ be a Noetherian $\FF_p$-algebra, and let $I$ denote the kernel of the multiplication map $R\otimes_{\FF_p}R\to R$. By Fogarty's result, we know that $R$ is $F$-finite, if and only if $\Omega_{R/\FF_p}=I/I^2$ is finitely generated. In this case, even though $I$ itself may not be finitely generated, it also satisfies a property that is typically proven using a finite set of generators: the families of ideals $\{I^n\}_{n\ge0}$ and $\{I^{[p^e]}\}_{e\ge0}$ are cofinal. This further highlights that morphisms of Noetherian $F$-finite rings share structural similarities with morphisms of finite type. For a more general statement, see \citep[Lemma 2.4]{QG24}.
\end{rem}

\section{The factorization for regular rings}
In this section, we show the existence of an $L$-smooth factorization for maps between regular $F$-finite rings. 

The following lemma plays a crucial role in what follows. First, it shows that $J$ is a regular ideal in the sense of Quillen \citep[]{Qui70}; see Remark \ref{regularideal}. Furthermore, it describes the associated graded pieces of the derived adic filtration of the Adams completion $\comp(R \otimes_{\mathbb{F}_p} S \to S)$ when $R$ is regular; see Lemma \ref{gradedpieces}. Finally, the finiteness of $J/J^2$ as an $S$-module is used to show that the completion of $R \otimes_{\mathbb{F}_p} S$ along the kernel of the map to $S$ is Noetherian and $F$-finite; see Proposition \ref{cplNoetherian}. 
\begin{lem}
\label{lemfgp}
Let $g:R\to S$ be a map of Noetherian $F$-finite rings, and let $J$ denote the kernel of the induced map $\nu:R\otimes_{\FF_p}S\to S$. If $R$ is regular, then
\[
L_{S/R\otimes_{\FF_p}S}[-1] \simeq Lg^*(L_{R/\FF_p})\simeq J/J^2[0]
\]
in $D(S)$, and $J/J^2$ is a finitely generated projective $S$-module.
\end{lem}
\begin{proof}
The factorization $S\to R\otimes_{\FF_p}S \to S$ of $\id_S$ yields the equivalence $L_{R\otimes_{\FF_p}S/S}\otimes^L_{R\otimes_{\FF_p}S}S\simeq L_{S/R\otimes_{\FF_p}S}[-1]$. Furthermore, the pushout square 
\[
\begin{tikzcd}
{\FF_p} \arrow[r] \arrow[d] & {S} \arrow[d] \\
{R} \arrow[r]           & {R\otimes_{\FF_p}S}          
\end{tikzcd}
\]
provides the equivalence $L_{R\otimes_{\FF_p}S/S}\simeq L_{R/\FF_p}\otimes^L_R R\otimes_{\FF_p}S$. Combining these, we obtain the first equivalence as stated in the lemma.
By Corollary \ref{regularislsmooth}, $R$ is $L$-smooth over $\FF_p$, and consequently, $Lg^*(L_{R/\FF_p})\simeq g^*(\Omega_{R/\FF_p})[0]$ is a finitely generated projective $S$-module. Thus, it suffices to show that $g^*(\Omega_{R/\FF_p})$ and $J/J^2$ are isomorphic. Let $I$ be the kernel of the multiplication map $\mu:R\otimes_{\FF_p}R\to R$. Consider the following commutative diagram:
\[
\begin{tikzcd}
{0} \arrow[r] & {I} \arrow[r] \arrow[d] & {R\otimes_{\FF_p}R} \arrow[r,"\mu"] \arrow[d, "\id\otimes g"] & {R} \arrow[r] \arrow[d, "g"] & {0} \\
{0} \arrow[r] & {J} \arrow[r]           & {R\otimes_{\FF_p}S} \arrow[r, "\nu"]           & {S} \arrow[r]           & {0}
\end{tikzcd}
\]
By the $2$-functoriality of $(-)^*$, we obtain the following isomorphisms:
\[
g^*(\Omega_{R/\FF_p})\cong g^*(\mu^*(I))\cong \nu^*((\id\otimes g)^*(I))
\]
Under the identification $R\otimes_{\FF_p} S \cong R\otimes_{\FF_p}R\otimes_R S$, the morphism $\id\otimes g$ corresponds to the map to the first factor $R\otimes_{\FF_p}R$. Hence, we have 
\[
(\id\otimes g)^*(I)\cong I \otimes_{R\otimes_{\FF_p} R} R\otimes_{\FF_p} R \otimes_R S\cong I\otimes_R S.
\]
Via the second factor inclusion, the above forms a commutative diagram of $R$-modules. By adjunction, this induces the following commutative diagram of $S$-modules:
\[
\begin{tikzcd}
{0} \arrow[r] & {I\otimes_R S} \arrow[r] \arrow[d] & {R\otimes_{\FF_p}S} \arrow[r] \arrow[d, "\id"] & {S} \arrow[r] \arrow[d, "\id"] & {0} \\
{0} \arrow[r] & {J} \arrow[r]           & {R\otimes_{\FF_p}S} \arrow[r]           & {S} \arrow[r]           & {0}
\end{tikzcd}
\]
Since $\Tor^R_1(R, S)$ vanishes, the first row of the diagram is exact, which implies that the first vertical map is an isomorphism. Therefore, we have $g^*(\Omega_{R/\FF_p})\cong \nu^*(J) \cong J/J^2$, as desired.
\end{proof}

\begin{rem}
The first equivalence in the previous lemma holds without any assumptions on $R$ and $S$, and can similarly be proven for arbitrary animated rings.
\end{rem}

\begin{rem}
\label{regularideal}
By \citep[Corollary 6.14]{Qui70}, the ideal $J$ in Lemma \ref{lemfgp} is a regular ideal in the sense of \citep[Definition 6.10]{Qui70}. This means that $J/J^2$ is a projective $S$-module and the canonical morphism of anticommutative algebras
\[
\Lambda^S_*(J/J^2)\to \Tor^{R\otimes_{\FF_p}S}_*(S,S)
\]
is an isomorphism.
\end{rem}

In \citep[Example 3.3.5]{BBST}, the following result is proven in the special case where $R=S$ and $g$ the identity map.

\begin{prop}
\label{cplNoetherian}
Let $g:R\to S$ be a map of Noetherian $F$-finite rings, and let $J$ denote the kernel of the induced map $R\otimes_{\FF_p}S\to S$. If $R$ is regular, then the $J$-adic completion $(R\otimes_{\FF_p} S)^{\wedge}_J$ is Noetherian and $F$-finite. Moreover, if $S$ is regular, then the completion is regular as well.
\end{prop}
\begin{proof}
We apply \citep[Corollaire 19.5.4]{EGAIV} with $A=S$, $B=R\otimes_{\FF_p}S$ and $C=S$. If $R$ is regular, $S\to R\otimes_{\FF_p} S$ is formally smooth by Proposition \ref{regularislsmooth}, as formal smoothness is preserved under base change. Consequently $(R\otimes_{\FF_p} S)^{\wedge}_J$ is naturally isomorphic to the adic completion of $\sym_S^{\bullet}(J/J^2)$ along the ideal of positive degree elements. It therefore suffices to show that $\sym_S^{\bullet}(J/J^2)$ is Noetherian and $F$-finite (and regular if $S$ is regular). 
Since $J/J^2$ is a finitely generated projective $S$-module, $\sym_S^{\bullet}(J/J^2)$ is a retract of a polynomial ring over $S$ in finitely many variables. In particular, $\sym_S^{\bullet}(J/J^2)$ is Noetherian and $F$-finite. If $S$ is regular, then the symmetric algebra, being a retract of a regular ring, is also regular by \citep[Corollary 1.11]{Cos77}.
\end{proof}

\begin{rem}
The proof of the previous proposition also shows that $(R\otimes_{\FF_p} S)^{\wedge}_J$ is complete along the kernel of the natural map to $S$.
\end{rem}

The next proposition shows that the natural map $R\otimes_{\FF_p} S \to (R\otimes_{\FF_p} S)^{\wedge}_J$ is both $L$-smooth and formally étale ("$L$-etale"). The proof requires the following lemma.

\begin{lem}
Let $R\to S$ be a map of Noetherian $F$-finite rings. Let $J$ denote the kernel of the induced map $R\otimes_{\FF_p}S\to S$, and $T=(R\otimes_{\FF_p} S)^{\wedge}_J$ be the $J$-adic completion. If $R$ is regular, then $L_{T/R\otimes_{\FF_p}S}$ is pseudo-coherent.
\end{lem}
\begin{proof}
Since both $R$ and $S$ are Noetherian $F$-finite rings, the cotangent complex
\[
L_{R\otimes_{\FF_p}S/\FF_p}\simeq (L_{R/\FF_p}\otimes_{R} R\otimes_{\FF_p}S) \oplus (L_{S/\FF_p}\otimes_{S} R\otimes_{\FF_p}S)
\]
is pseudo-coherent. If $R$ is regular, then by Proposition \ref{cplNoetherian}, $T$ is Noetherian and $F$-finite, which implies that $L_{T/\FF_p}$ is also pseudo-coherent. The result then follows from the transitivity fiber sequence induced by the maps $\FF_p\to R\otimes_{\FF_p}S \to T$, together with the “two-out-of-three” property of pseudo-coherence in a fiber sequence.
\end{proof}

\begin{prop}
\label{vanishingcotangent}
Let $g:R\to S$ be a map of Noetherian $F$-finite rings, and let $J$ denote the kernel of the induced map $R\otimes_{\FF_p}S\to S$. If $R$ is regular, then the cotangent complex of the natural map
\[
R\otimes_{\FF_p} S \to (R\otimes_{\FF_p} S)^{\wedge}_J=T
\]
is acyclic.
\end{prop}
\begin{proof}
Assume that $R$ is regular, and let $I$ denote the kernel of the natural map $T\to S$. Since $T$ is $I$-adically complete and $L_{T/R\otimes_{\FF_p}S}$ is pseudo-coherent in $D(T)$, the latter is derived $I$-complete by \cite[\href{https://stacks.math.columbia.edu/tag/0A05}{Lemma 0A05}]{stacks-project}. Consequently, \cite[\href{https://stacks.math.columbia.edu/tag/0922}{Proposition 0922}]{stacks-project} provides the equivalence
\[
L_{T/R\otimes_{\FF_p}S}\simeq \lim_n L_{T/R\otimes_{\FF_p}S}\otimes^L_T T/I^n.
\]
Since $T$ is Noetherian, the ideal $I$ is pro Tor-unital by \citep[Theorem 2.3]{Mor18}. Thus, by \citep[Remark 4.6 (ii)]{Mor18}, the following equivalence holds:
\[
\lim_n L_{T/R\otimes_{\FF_p}S}\otimes_T T/I^n \simeq \lim_n L_{T/I^n /R\otimes_{\FF_p}S} 
\]
Let $I_+=\ker(\sym_S^{\bullet}(J/J^2)\to S)$ denote the ideal of positive degree elements. Then
\[
T/I^n\cong\sym_S^{\bullet}(J/J^2)/I_+^n\cong (R\otimes_{\FF_p}S)/J^n,
\]
where the second isomorphism is shown in the proof of \citep[Corollaire 19.5.4]{EGAIV} on pages 93 and 94. Consequently, 
\[
\lim_n L_{T/I^n /R\otimes_{\FF_p}S}\simeq \lim_n L_{(R\otimes_{\FF_p}S)/J^n /R\otimes_{\FF_p}S}.
\]
According to Lemma \ref{lemfgp} and \citep[Example 1.4]{Mor18}, the ideal $J$ is also pro Tor-unital. Thus, by \citep[Remark 4.6 (i)]{Mor18}, the limit $\lim_n L_{(R\otimes_{\FF_p}S)/J^n /R\otimes_{\FF_p}S}$ is acyclic, completing the proof.
\end{proof}

\begin{cor}
\label{regfact}
Let $g:R\to S$ be a map of regular Noetherian $F$-finite $\FF_p$-algebras. Let $J$ denote the kernel of the induced map $R\otimes_{\FF_p}S\to S$, and let $T=(R\otimes_{\FF_p} S)^{\wedge}_J$ be the $J$-adic completion. Then the first-factor inclusion $R\to R\otimes_{\FF_p} S$ induces a factorization
\[
\begin{tikzcd}
                         & {T} \arrow[rd, "\psi"] &    \\
{R} \arrow[ru, "\varphi"] \arrow[rr, "g"] &               & {S}
\end{tikzcd}
\]
where $\varphi$ is $L$-smooth with cotangent complex $L_{T/R}\simeq L_{S/\FF_p}\otimes^L_{S}T$ and $\psi$ is surjective.
\end{cor}
\begin{proof}
Since $S$ is regular and $F$-finite, it is $L$-smooth over $\FF_p$ by Corollary \ref{regularislsmooth}. Consequently, by Proposition \ref{basechangelsmooth}, the base change $R\to R\otimes_{\FF_p} S$ is also $L$-smooth. Now, as the cotangent complex of the map $R\otimes_{\FF_p}S\to T$ vanishes by Proposition \ref{vanishingcotangent}, the canonical fiber sequence induced by $R\to R\otimes_{\FF_p}S\to T$ gives
\[
L_{T/R}\simeq L_{(R\otimes_{\FF_p}S)/R}\otimes^L_{R\otimes_{\FF_p}S}T \simeq L_{S/\FF_p}\otimes^L_{S}T.
\]
In particular, $\varphi$ is $L$-smooth by Corollary \ref{regularislsmooth}.
\end{proof}

\begin{rem}
Under the same assumptions as in the previous corollary, we can provide an alternative proof that $T$ is regular: Since $\FF_p\to R$ and $R\to T$ are both $L$-smooth, the same is true for their composition $\FF_p\to T$. As $T$ is Noetherian and $F$-finite, it follows that $T$ is regular by Corollary \ref{regularislsmooth}.
\end{rem}

As surjective maps of regular rings are regular immersions, we immediately obtain the following: 

\begin{cor}
Any map of regular Noetherian $F$-finite rings $R\to S$ can be factored into an $L$-smooth map $R\to T$ of regular Noetherian $F$-finite rings and a regular immersion $T\to S$. Any such factorization induces an exact sequence
\[
0\to H^{-1}(L_{S/R})\to I/I^2 \to \Omega_{T/R}\otimes_T S \to \Omega_{S/R}\to 0,
\]
which describes all the nonvanishing cohomology groups of $L_{S/R}$.
\end{cor}
\begin{proof}
By Corollary \ref{regfact}, there exists a factorization $R\to T\to S$ with the desired properties. The kernel $I$ of the second map is regular by \cite[\href{https://stacks.math.columbia.edu/tag/0E9J}{Lemma 0E9J}]{stacks-project}. In particular, this gives a canonical equivalence $L_{S/T}\simeq I/I^2[1]$ \cite[\href{https://stacks.math.columbia.edu/tag/08SJ}{Lemma 08SJ}]{stacks-project}. 

Since $R$ and $S$ are $L$-smooth over $\FF_p$ by Corollary \ref{regularislsmooth}, the canonical fiber sequence of cotangent complexes associated with $\FF_p\to R\to S$ shows that $L_{S/R}$ is a perfect complex concentrated in cohomological degrees $0$ and $-1$. Therefore, the long exact sequence associated with the fiber sequence of cotangent complexes induced by the factorization $R\to T\to S$ is
\[
0\to H^{-1}(L_{S/R})\to I/I^2 \to \Omega_{T/R}\otimes_T S \to \Omega_{S/R}\to 0.
\]
\end{proof}

Although not required for the subsequent discussion, we present the following statement as its proof relies on similar results from \citep{Mor18} and uses arguments analogous to those in the previous proposition.

\begin{prop}
Let $R$ be a Noetherian $F$-finite ring and $I\subset R$ an ideal. Then, the cotangent complex of the canonical map $R\to R^{\wedge}_I$ is acyclic.
\end{prop}
\begin{proof}
Since $R$ and $R^{\wedge}_I$ are both Noetherian and $F$-finite, the cotangent complex $L_{R^{\wedge}_I/R}$ is pseudo-coherent. Let $J$ denote the kernel of the natural map $R^{\wedge}_I\to R/I$. As $R^{\wedge}_I$ is $J$-complete, we obtain an equivalence
\[
L_{R^{\wedge}_I/R}\simeq \lim_n L_{R^{\wedge}_I/R}\otimes^L_{R^{\wedge}_I}R^{\wedge}_I/J^n.
\]
Since $R^{\wedge}_I$ is Noetherian, $J$ is pro Tor-unital by \citep[Theorem 2.3]{Mor18}. Consequently 
\[
\lim_n L_{R^{\wedge}_I/R}\otimes^L_{R^{\wedge}_I}R^{\wedge}_I/J^n \simeq \lim_n L_{(R^{\wedge}_I/J^n)/R} \simeq \lim_n L_{(R/I^n)/R}.
\]
Finally, the last expression is acyclic because $I$ is pro Tor-unital as an ideal of a Noetherian ring.
\end{proof}

\section{The factorization in general}
The following result plays a key role in reducing the general case to the regular case. The proof provides an explicit construction of the pushout square in question.

\begin{prop}
\label{pushoutrings}
Any map $f:R\to S$ of Noetherian $F$-finite $\FF_p$-algebras can be extended to a pushout square
\[
\begin{tikzcd}
{R'} \arrow[r, "g"] \arrow[d] & {S'} \arrow[d] \\
{R} \arrow[r, "f"]           & {S}          
\end{tikzcd}
\]
of $\FF_p$-algebras such that: 
\begin{enumerate}
\item the vertical maps are surjective,
\item $R'$ and $S'$ are regular Noetherian and $F$-finite, and
\item $R'$ is complete along the ideal $\ker(R'\to R)$, and $S'$ is complete along $\ker(S'\to S)$.
\end{enumerate}
\end{prop}
\begin{proof}
We choose a finite family of $p$-generators for $R$ and extend their images under $f$ to form a finite family of $p$-generators for $S$. In this way, $f$ can be lifted at each step in the construction of the regular $F$-finite rings, as described in \citep[Remark 13.6]{Gab04}. Consequently, there exists a commutative diagram
\[
\begin{tikzcd}
{G(R)} \arrow[r, "G(f)"] \arrow[d] & {G(S)} \arrow[d] \\
{R} \arrow[r, "f"] & {S}
\end{tikzcd}
\]
in which the vertical maps are surjective, and $G(R)$ and $G(S)$ are regular Noetherian and $F$-finite; see \cite[Proposition 2.3.3]{BBST}. However, in general, $\ker(G(R)\to R)G(S)\subset \ker(G(S)\to S)$ might be a strict inclusion, in which case this diagram is not cocartesian. Let $\{h_1,\dots,h_n\}$ be a finite set of generators for $\ker(G(S)\to S)$, and let $\tilde{R}=G(R)[X_1,\dots,X_n]$ be a polynomial ring in $n$ variables over $G(R)$. Then $G(f)$ induces unique maps $\tilde{R}\to G(S)$ and $\tilde{R}\to R$, such that the polynomial variable $X_i$ is mapped to $h_i$ in $G(S)$ and to $0$ in $R$ for $i=1,\dots,n$. Now we have
\[
J=\ker(\tilde{R}\to R)=\ker(G(R)\to R)\tilde{R} + (x_1,\dots,x_n)
\]
and consequently, the desired equality $JG(S)=\ker(G(S)\to S)$ holds. Thus, the induced diagram 
\[
\begin{tikzcd}
{\tilde{R}} \arrow[r] \arrow[d] & {G(S)} \arrow[d] \\
{R} \arrow[r, "f"] & {S}
\end{tikzcd}
\]
is a pushout square that satisfies the properties (1) and (2) stated in the proposition. Moreover, by \citep[Construction 2.2.3]{BBST}, $G(S)$ is complete along $JG(S)=\ker(G(S)\to S)$. Consequently, the maps $\tilde{R}\to G(S)$ and $\tilde{R}\to R$ factor uniquely over $\tilde{R}\to\tilde{R}^{\wedge}_J$, and the induced diagram satisfies all the desired properties (1)-(3).
\end{proof}

\begin{rem}
The pushout square constructed in the proof of Proposition \ref{pushoutrings} is, in general, not a pushout square of animated rings, as $R$ and $S'$ need not be Tor-independent over $R'$.
\end{rem}

\begin{con}
\label{classicconstruction}
Let $f:R\to S$ be a map of Noetherian $F$-finite $\FF_p$-algebras. Choose a pushout square of $\FF_p$-algebras 
\[
\begin{tikzcd}
{R'} \arrow[r, "g"] \arrow[d] & {S'} \arrow[d] \\
{R} \arrow[r, "f"]           & {S}          
\end{tikzcd}
\]
satisfying properties (1) and (2) of Proposition \ref{pushoutrings}. Let $J'$ denote the kernel of the map $R'\otimes_{\FF_p}S'\to S'$ induced by $g$ and the identity on $S'$, and denote $T'=(R'\otimes_{\FF_p}S')^{\wedge}_{J'}$ the $J'$-adic completion. This construction yields a commutative diagram:
\[
\begin{tikzcd}
                         & {T'} \arrow[rd, "\psi{'}"] &    \\
{R'} \arrow[ru, "\varphi{'}"] \arrow[rr, "g"] &               & {S'}
\end{tikzcd}
\]
where $\varphi'$ is induced by the first factor inclusion, and $\psi'$ is the canonical projection. Changing the base along $R'\to R$ induces a commutative diagram 
\[
\begin{tikzcd}
                         & {T} \arrow[rd, "\psi"] &    \\
{R} \arrow[ru, "\varphi"] \arrow[rr, "f"] &               & {S.}
\end{tikzcd}
\]
\end{con}

\begin{thm}
\label{maintheorem}
In the context of Construction \ref{classicconstruction}, the $\FF_p$-algebra $T$ is Noetherian and $F$-finite, $\varphi$ is $L$-smooth, and $\psi$ is surjective.
\end{thm}
\begin{proof}
By Proposition \ref{cplNoetherian}, $T'$ is Noetherian and $F$-finite. Since the induced map $T'\to T$ is surjective, it follows that $T$ is also Noetherian and $F$-finite. By Corollary \ref{regfact}, $\varphi'$ is $L$-smooth, and the natural projection $\psi'$ is surjective. As $\varphi'$ is an $L$-smooth morphism of Noetherian rings, it is flat. This implies that $T'$ and $R$ are Tor-independent over $R'$. Consequently, by Proposition \ref{basechangelsmooth}, the map $\varphi:R\to R\otimes_{R'} T'=T$ is also $L$-smooth. Finally, since $\psi'$ is surjective, the map $\psi:T\to[\id_R\otimes\psi'] R\otimes_{R'}S'\cong S$ is also surjective.
\end{proof}

\begin{rem}
With the same notation as in Construction \ref{classicconstruction} and the proof of Theorem \ref{maintheorem}, the cotangent complex of $\varphi :R\to T$ is given by 
\[
L_{T/R}\simeq L_{T'/R'}\otimes^L_{T'} T \simeq L_{S'/\FF_p}\otimes^L_{S'}T.
\]
Here, the first equivalence holds by \cite[\href{https://stacks.math.columbia.edu/tag/08QQ}{Lemma 08QQ}]{stacks-project}, and the second equivalence is a consequence of Corollary \ref{regfact}.
\end{rem}

\section{A Description via the Adams completion}
In this section, we use the formalism of Adams completion, as developed in the articles \citep{Car08}, \citep{bha12}, and \citep{BBST}, to provide an alternative description of $T$ as stated in Theorem \ref{maintheorem}. Many of the arguments draw inspiration from \citep[Sections 3.2 and 3.3]{BBST}.

Let $A \to B$ be a map of animated rings with connective fiber $I$. As shown in \citep[Proposition 3.2.6]{BBST}, the \emph{Adams completion} of this map
\[
\comp(A\to B)=\lim_n A/I^{\otimes^L_A n}
\]
is naturally equipped with the so-called \emph{derived $I$-adic filtration}. This is a complete descending $\NN$-indexed filtration with associated graded pieces 
\[
\gr^n\comp(A\to B)\simeq L\sym^n_B(L_{B/A}[-1])
\]
for $n\ge0$.

\begin{lem}
\label{gradedpieces}
Let $g:R\to S$ be a map of Noetherian $F$-finite rings, and let $J$ denote the kernel of the induced map $R\otimes_{\FF_p}S\to S$. If $R$ is regular, then there is a natural identification
\[
\gr^n\comp(R\otimes_{\FF_p}S\to S)\simeq J^n/J^{n+1}[0]
\]
for all $n\ge0$. 
\end{lem}
\begin{proof}
By Lemma \ref{lemfgp}, there is a natural equivalence $L_{S/R\otimes_{\FF_p}S}[-1]\simeq J/J^2[0]$, and $J/J^2$ is projective. This induces an equivalence
\[
L\sym^n_S(L_{S/R\otimes_{\FF_p}S}[-1])\simeq \sym^n_S(J/J^2)[0]\simeq J^n/J^{n+1}[0]
\]
by \citep[Corollaire 19.5.4]{EGAIV} applied to $A=S, B=R\otimes_{\FF_p}S$ and $C=S$.
\end{proof}

We use the identification of the associated graded pieces of the derived filtration to prove that the Adams completion is isomorphic to the classical completion. This was shown in \citep[Example 3.3.5]{BBST} for the special case $R=S$ and $g$ being the identity map.
\begin{prop}
\label{adamsequalsclassical}
Let $g:R\to S$ be a map of Noetherian $F$-finite rings, and let $J$ denote the kernel of the induced map $\nu:R\otimes_{\FF_p}S\to S$. If $R$ is regular, then the canonical map
\[
\comp(R\otimes_{\FF_p}S\to S)\to (R\otimes_{\FF_p}S)^{\wedge}_J
\]
is an equivalence. In particular, $\comp(R\otimes_{\FF_p}S\to S)$ is discrete.
\end{prop}
\begin{proof}
The equivalence $\comp(R\otimes_{\FF_p}S\to S)\simeq \lim_n((R\otimes_{\FF_p}S)/J^{\otimes n})$ induces a surjective map 
\[
\pi_0(\comp(R\otimes_{\FF_p}S\to S)) \to \lim_n \pi_0((R\otimes_{\FF_p}S)/J^{\otimes n}),
\]
and by analyzing the associated long exact homotopy sequences, the target can be identified with $(R\otimes_{\FF_p}S)^{\wedge}_J$. This yields a canonical map
\[
\comp(R\otimes_{\FF_p}S\to S)\to\pi_0(\comp(R\otimes_{\FF_p}S\to S))\to (R\otimes_{\FF_p}S)^{\wedge}_J.
\]
This morphism is filtered with respect to the complete derived adic filtration on $\comp(R\otimes_{\FF_p}S\to S)$ and the classical complete descending filtration on $(R\otimes_{\FF_p}S)^{\wedge}_J$. To show that it is an equivalence, it therefore suffices to check that the induced maps on associated graded pieces are equivalences. But by Lemma \ref{gradedpieces}, the graded pieces of both filtrations are naturally identified.
\end{proof}

\begin{cor}
\label{regadamsnoetherianffinite}
Let $R\to S$ be a map of Noetherian $F$-finite $\FF_p$-algebras. If $R$ is regular, then the Adams completion $\comp(R\otimes_{\FF_p}S\to S)$ is a Noetherian $F$-finite ring. Furthermore, if both $R$ and $S$ are regular, the Adams completion is regular as well.  
\end{cor}
\begin{proof}
This follows directly from Proposition \ref{cplNoetherian} and Proposition \ref{adamsequalsclassical}.
\end{proof}

\begin{lem}
\label{pushoutcompletions}
Let 
\[
\begin{tikzcd}
{R'} \arrow[r, "g"] \arrow[d] & {S'} \arrow[d] \\
{R} \arrow[r, "f"]           & {S}          
\end{tikzcd}
\]
be a pushout square of Noetherian $F$-finite $\FF_p$-algebras in which the vertical maps are finite. If $R'$ is regular, then 
\[
\begin{tikzcd}
{R'} \arrow[r] \arrow[d] & {\comp(R'\otimes_{\FF_p}S'\to S')} \arrow[d] \\
{R} \arrow[r]           & {\comp(R\otimes_{\FF_p}S'\to S)}          
\end{tikzcd}
\]
is a pushout square of animated rings.
\end{lem}
\begin{proof}
If $R'$ is regular, then $R\in D(R')$ is a perfect complex. This implies the natural equivalence
\[
R\otimes^L_{R'}\comp(R'\otimes_{\FF_p}S'\to S')\simeq \comp(R\otimes_{\FF_p} S'\to R\otimes^L_{R'} S')
\]
of animated rings. By assumption, $\pi_0(R\otimes^L_{R'} S')=R\otimes_{R'} S'$ is canonically isomorphic to $S$. Therefore, we have 
\[
\comp(R\otimes_{\FF_p} S'\to R\otimes^L_{R'} S')\simeq \comp(R\otimes_{\FF_p} S'\to S)
\]
by \citep[Proposition 3.2.8 and 3.2.10]{BBST}. 
\end{proof}

\begin{rem}
If the initial diagram in the previous lemma is not cocartesian, there still exists a natural equivalence of animated rings:
\[
R\otimes^L_{R'}\comp(R'\otimes_{\FF_p}S'\to S)\simeq \comp(R \otimes_{\FF_p}S'\to S)
\]
\end{rem}

\begin{con}
\label{derivedconstruction}
Let $f:R\to S$ be a map of Noetherian $F$-finite $\FF_p$-algebras. Begin by selecting a finite family of $p$-generators for $R$, and extend their images under $f$ to a finite family of $p$-generators for $S$. Using these $p$-generators, construct the regular Noetherian and $F$-finite ring $S'=G(S)$ together with the surjective map $q: S' \to S$, as described in \citep[Remark 13.6]{Gab04} and \citep{BBST}.
Define $T=\comp(R \otimes_{\FF_p} S' \to S)$ as the Adams completion of the map induced by $f$ and $q$. This yields a natural factorization of animated rings
\[
\begin{tikzcd}
                         & {T} \arrow[rd, "\psi"] &    \\
{R} \arrow[ru, "\varphi"] \arrow[rr, "f"] &               & {S.}
\end{tikzcd}
\]
\end{con}

\begin{thm}
The factorization of $f$ in Construction \ref{derivedconstruction} is an $L$-smooth-by-surjective factorization, which means that $T$ is a discrete Noetherian $F$-finite ring, the map $\varphi$ is $L$-smooth, and $\psi$ is surjective.
\end{thm}
\begin{proof}
Starting with a finite family of $p$-generators for $R$ in Construction \ref{derivedconstruction} ensures that there exists a pushout square of $\FF_p$-algebras
\[
\begin{tikzcd}
{R'} \arrow[r, "g"] \arrow[d] & {S'} \arrow[d] \\
{R} \arrow[r, "f"]           & {S}          
\end{tikzcd}
\]
that satisfies properties (1) and (2) of Proposition \ref{pushoutrings}, as explained in its proof.
By Proposition \ref{adamsequalsclassical}, the canonical map
\[
\comp(R'\otimes_{\FF_p}S'\to S')\to (R'\otimes_{\FF_p}S')^{\wedge}_{J'}=T'
\]
is an equivalence, where $J'$ denotes the kernel of $R'\otimes_{\FF_p}S'\to S'$. Applying Lemma \ref{pushoutcompletions}, we deduce that the diagram 
\[
\begin{tikzcd}
{R'} \arrow[r] \arrow[d] & {T'} \arrow[d] \\
{R} \arrow[r]           & {\comp(R\otimes_{\FF_p}S'\to S)}          
\end{tikzcd}
\]
is a pushout square of animated rings. By Corollary \ref{regfact}, the map $R'\to T'$ is $L$-smooth and therefore flat, which implies that the diagram is also a pushout square of $\FF_p$-algebras. Thus, $\comp(R\otimes_{\FF_p}S'\to S)$ is discrete, the canonical map $\varphi$ is $L$-smooth, and $\psi$ is surjective as explained in the proof of Theorem \ref{maintheorem}.
\end{proof}

\begin{cor}
The factorization constructed in Construction \ref{classicconstruction} is independent of the choice of $R'$.
\end{cor}
\begin{proof}
This follows directly from the fact that Construction \ref{derivedconstruction} is independent of $R'$.
\end{proof}

\begin{cor}
Let $f:R\to S$ be a map of Noetherian $F$-finite rings, and 
\[
\begin{tikzcd}
{R'} \arrow[r] \arrow[d] & {S'} \arrow[d] \\
{R} \arrow[r]           & {S}          
\end{tikzcd}
\]
a pushout square satisfying conditions (1) and (2) of Proposition \ref{pushoutrings}. Denote by $J$ the kernel of the induced map $R\otimes_{\FF_p}S'\to S$. Then the $J$-adic completion $(R\otimes_{\FF_p}S')^{\wedge}_J$ is Noetherian and $F$-finite.
\end{cor}
\begin{proof}
The Milnor short exact sequence on $\pi_0$ 
shows that $(R\otimes_{\FF_p}S')^{\wedge}_J$ is a quotient of the (discrete) Noetherian and $F$-finite ring $\comp(R\otimes_{\FF_p}S'\to S)$.
\end{proof}

The next result describes the cotangent complex of the second map in the $L$-smooth factorization.

\begin{prop}
\label{cotangentcompletion}
Let 
\[
\begin{tikzcd}
{R'} \arrow[d] \arrow[r] & {S'} \arrow[d] \\
{R} \arrow[r]           & {S}          
\end{tikzcd}
\]
be a pushout square of Noetherian $F$-finite $\FF_p$-algebras, where $R'$ and $S'$ are regular, and the vertical maps are surjective. Then the natural map
\[
L_{S/R\otimes_{\FF_p}S'}\to L_{S/\comp(R\otimes_{\FF_p}S'\to S)}
\]
is an equivalence.
\end{prop}
\begin{proof}
The sequence of maps $R\otimes_{\FF_p}S'\to\comp(R\otimes_{\FF_p}S'\to S)\to S$ induces a fiber sequence
\[
L_{\comp(R\otimes_{\FF_p}S'\to S)/R\otimes_{\FF_p}S'}\otimes^L_{\comp(R\otimes_{\FF_p}S'\to S)}S
\to
L_{S/R\otimes_{\FF_p}S'}
\to
L_{S/\comp(R\otimes_{\FF_p}S'\to S)}.
\]
Thus, it suffices to show that $L_{\comp(R\otimes_{\FF_p}S'\to S)/R\otimes_{\FF_p}S'}\otimes^L_{\comp(R\otimes_{\FF_p}S'\to S)}S$ is acyclic.
By Lemma \ref{pushoutcompletions} and the “two-out-of-three” property for cocartesian squares, the diagram
\[
\begin{tikzcd}
{R'\otimes_{\FF_p} S'} \arrow[r] \arrow[d] & {\comp(R'\otimes_{\FF_p}S'\to S')} \arrow[d] \\
{R\otimes_{\FF_p} S'} \arrow[r]           & {\comp(R\otimes_{\FF_p}S'\to S)}          
\end{tikzcd}
\]
is a pushout square of animated rings. Consequently, we have an equivalence
\[
L_{\comp(R\otimes_{\FF_p}S'\to S)/R\otimes_{\FF_p}S'}\otimes^L_{\comp(R\otimes_{\FF_p}S'\to S)}S
\simeq 
L_{\comp(R'\otimes_{\FF_p}S'\to S)/R'\otimes_{\FF_p}S'}\otimes^L_{\comp(R'\otimes_{\FF_p}S'\to S)}S.
\]
Since the map $\comp(R'\otimes_{\FF_p}S'\to S)\to S$ factors through $S'$, it suffices to show that
\[
L_{\comp(R'\otimes_{\FF_p}S'\to S)/R'\otimes_{\FF_p}S'}\otimes^L_{\comp(R'\otimes_{\FF_p}S'\to S)}S'
\]
is acyclic. Hence, we may assume $R=R'$ and $S=S'$ are regular Noetherian and $F$-finite rings. 
By Corollary \ref{regadamsnoetherianffinite}, $\comp(R\otimes_{\FF_p}S\to S)$ is a regular $F$-finite ring. This implies that the kernel $I$ of the map
$\comp(R\otimes_{\FF_p}S\to S)\to S$
is a regular ideal. Consequently, by \citep[Corollary 6.14]{Qui70} (resp. \cite[\href{https://stacks.math.columbia.edu/tag/08SJ}{Lemma 08SJ}]{stacks-project}), we have an equivalence
\[
L_{S/\comp(R\otimes_{\FF_p}S\to S)}
\simeq
I/I^2[1].
\]
If $J$ denotes the kernel of $R\otimes_{\FF_p}S\to S$, we also have the equivalence
\[
L_{S/R\otimes_{\FF_p}S}
\simeq
J/J^2[1]
\]
by Lemma \ref{gradedpieces}. The proposition follows because $J/J^2$ and $I/I^2$ are canonically isomorphic to $I_+/I_+^2$, where $I_+=\ker(\sym_S^{\bullet}(J/J^2)\to S)$.
\end{proof}

\emergencystretch=1.1em
\printbibliography
\end{document}